\documentclass[11pt]{article}
\usepackage{amssymb}
\usepackage{amscd}
\usepackage{amsthm}
\usepackage{amsmath}
\usepackage{latexsym}
\usepackage{graphicx}
\usepackage{delarray}
\usepackage{epic}
\usepackage{float}

\begin{document}

%%%%%%%%%%%%%%%%%%%%%%%%%%%%%%%%%%%%%%%%%%%%%%%%%%%%%%%%%%%%%%%%%%%%%%%%

\newtheorem{theorem}{Theorem}[section]
\newtheorem{definition}[theorem]{Definition}
\newtheorem{lemma}[theorem]{Lemma}
\newtheorem{proposition}[theorem]{Proposition}
\newtheorem{corollary}[theorem]{Corollary}
\newtheorem{example}[theorem]{Example}
\newtheorem{remark}[theorem]{Remark}

\hfuzz5pt % Don't bother to report over-full boxes if over-edge is < 5pt

%%%%%%%%%%%%%%%%%%%%%%%%%%%%%%%%%%%%%%%%%%%%%%%%%%%%%%%%%%%%%%%%%%%%%%%%

%% macros: Tobias

\newcommand{\gt}{\tilde{g}}
\newcommand{\R}{\mathbb{R}}
\newcommand{\C}{\mathbb{C}}
\newcommand{\Z}{\mathbb{Z}}
\newcommand{\Zd}{\mathbb{Z}^d}
\newcommand{\Ztd}{\mathbb{Z}^{2d}}
\newcommand{\Rd}{\R^d}
\newcommand{\Rtd}{\R^{2d}}
\newcommand{\Zp}{Z_p}
\newcommand{\al}{\alpha}
\newcommand{\be}{\beta}
\newcommand{\om}{\omega}
\newcommand{\ga}{\gamma}
\newcommand{\la}{\lambda}
\newcommand{\La}{\Lambda}
\newcommand{\ala}{\la^\circ}
\newcommand{\aLa}{\La^\circ}
\newcommand{\nat}{\natural}
\newcommand{\G}{\mathcal{G}}
\newcommand{\A}{\mathcal{A}}
\newcommand{\M}{\mathcal{M}}
\newcommand{\ka}{\kappa}
\newcommand{\cast}{\circledast}
\newcommand{\id}{\mbox{Id}}

\newcommand{\conv}[2]{{#1}\,\ast\,{#2}}
\newcommand{\twc}[2]{{#1}\,\nat\,{#2}}
\newcommand{\mconv}[2]{{#1}\,\cast\,{#2}}
\newcommand{\set}[2]{\Big\{ \, #1 \, \Big| \, #2 \, \Big\}}
\newcommand{\inner}[2]{\langle #1,#2\rangle}
\newcommand{\dotp}[2]{ #1 \, \cdot \, #2}

%% macros: Ewa

\newcommand{\Zpd}{\Zp^d}
\newcommand{\I}{\mathcal{I}}
\newcommand{\Zq}{Z_q}
\newcommand{\Zqd}{Zq^d}
\newcommand{\Zak}{\mathcal{Z}_{\alpha}}
\newcommand{\D}{\mathcal{D}}

%%%%%%%%%%%%%%%%%%%%%%%%%%%%%%%%%%%%%%%%%%%%%%%%%%%%%%%%%%%%%%%%%%%%%%%%
%%%%%%%%%%%%%%%%%%%%%%%%%%%%%%%%%%%%%%%%%%%%%%%%%%%%%%%%%%%%%%%%%%%%%%%%

\title{A Constructive Inversion Framework for Twisted Convolution}

\author{Yonina C.~Eldar\footnote{Dept.~of Electrical Engineering,
Technion--Israel Institute of Technology, 32000 Haifa, Israel.
\newline Tel.: +972-4-8293256, Fax.: +972-4-8295757,
Email: \texttt{yonina@ee.technion.ac.il}}\,,
\, Ewa Matusiak\footnote{Faculty~of Math., University of Vienna, 1090
Vienna, Austria. Tel.: +43-1-4277-50693, Fax.: +43-1-4277-50690,
\newline Email: \texttt{ewa.matusiak@univie.ac.at},
\texttt{tobias.werther@univie.ac.at}. }\,\,, \, Tobias Werther $^\dagger$}

%% \date{}

\maketitle

%%%%%%%%%%%%%%%%%%%%%%%%%%%%%%%%%%%%%%%%%%%%%%%%%%%%%%%%%%%%%%%%%%%%%%%%


 \vspace{.2cm}

%%%%%%%%%%%%%%%%%%%%%%%%%%%%%%%%%%%%%%%%%%%%%%%%%%%%%%%%%%%%%%%%%%%%%%%%
%%%%%%%%%%%%%%%%%%%%%%%%%%%%%%%%%%%%%%%%%%%%%%%%%%%%%%%%%%%%%%%%%%%%%%%%

\begin{abstract}
In this paper we develop constructive invertibility conditions
for the twisted convolution. Our approach is based on splitting the
twisted convolution with rational parameters into a finite number of
weighted convolutions, which can be interpreted as another twisted
convolution on a finite cyclic group. In analogy with the twisted
convolution of finite discrete signals, we derive an anti-homomorphism
between the sequence space and a suitable matrix algebra which preserves
the algebraic structure. In this way, the problem reduces to the
analysis of finite matrices whose entries are sequences supported on
corresponding cosets. The invertibility condition then follows from
Cramer's rule and Wiener's lemma for this special class of matrices.
The problem results from a well known approach of studying
the invertibility properties of the Gabor frame operator in the rational
case. The presented approach gives further insights into Gabor frames.
In particular, it can be applied for both the continuous
(on $\Rd$) and the finite discrete setting. In the latter case,
we obtain algorithmic schemes for directly computing the inverse of Gabor
frame-type matrices equivalent to those known in the literature.
\end{abstract}

%%%%%%%%%%%%%%%%%%%%%%%%%%%%%%%%%%%%%%%%%%%%%%%%%%%%%%%%%%%%%%%%%%%%%%%%
%%%%%%%%%%%%%%%%%%%%%%%%%%%%%%%%%%%%%%%%%%%%%%%%%%%%%%%%%%%%%%%%%%%%%%%%

\section{Introduction}

%% In this section we give an outline of the manuscript and
%% explain why we treat the invertibility of a twisted convolution
%% operator.

Twisted convolution arises naturally in the context of time frequency
operators, more specifically in the treatment of Gabor frames
\cite{CH03,Gro01}. The study of inversion schemes of twisted convolution has,
therefore, a major impact on the analysis of Gabor frames.
Our method is originated by the Janssen representation
of Gabor frame operators \cite{Jan95} and simplifies the approach given in
\cite{WEN05}. A different, however, equivalent method for studying Gabor
frame operators is the well known Zibulski-Zeevi representation \cite{ZZ97}
based on a generalized Zak-transform.

In contrast to the standard convolution, the twisted convolution is
not commutative. This is opposed to the possibility of applying powerful
tools from harmonic analysis, such as Wiener's Lemma, in order to study
twisted convolution operators. Recently, in \cite{WEN05}, the authors
described an new approach to classify the invertibility of
$\ell^1$-sequences with respect to the twisted convolution for rational
parameters.

In this manuscript we extend the idea of \cite{WEN05} in the sense
that we take a different approach which allows far better insights
into the problem. Specifically, we only deal with sequences and
show explicitly how efficient inversion schemes can be derived by
rather simple (though sophisticated) manipulations of the twisted
convolution. The essential idea is to split up
the twisted convolution into a finite number of sums that can be
incorporated into a special matrix algebra. In this matrix algebra we then
prove a special type of Wiener Lemma which is the most challenging
part from a mathematical perspective.

The paper is organized as follows. The first section briefly outlines
the basic definition of the twisted convolution. In this section we
further discuss the example of twisted convolution on the finite group
$Z_p\times Z_p$. This example serves the purpose to motivate the
introduction of the matrix algebra that appears in Section 3
where we prove Wiener's Lemma for a special subalgebra. Section 4
links the twisted convolution to time-frequency operators. More
specifically, it shows how the results shown in Section 3 can
be used in the context of Gabor frames. In the last section, we
give a short outline of the application of the presented approach
for inverting frame-like Gabor operators.

%%%%%%%%%%%%%%%%%%%%%%%%%%%%%%%%%%%%%%%%%%%%%%%%%%%%%%%%%%%%%%%%%%%%%%%%

\section{Twisted Convolution}

%% We define the twisted convolution and emphasize the main
%% difficulty which arises from the non-commutativity due
%% to the additional phase factor (in contrast to the standard
%% convolution) and time-frequency operators!!!

Let $p$ and $q$ be integers and relatively prime. We define the {\it twisted
convolution} for sequences $a,b\in\ell^1(\Ztd)$ by
\begin{equation}\label{eq:twconv}
(\twc{a}{b})_{m,n} = \sum_{k,l\in\Zd} a_{k,l}b_{m-k,n-l}\om^{\dotp{(m-k)}{l}}
\end{equation}
where $\om = e^{2\pi i q/p}$ and $\cdot$ denotes the inner product
in $\Rd$. Although the twisted convolution
depends on $p,q$ we do not specify this dependence because $p,q$
will always be given and fixed beforehand. In Section 4 we show
how the twisted convolution is related to a class of operators
with a special time-frequency representation.

In contrast to the conventional convolution with symbol $\ast$,
in which $\om=1$, the twisted convolution is not commutative, and turns
$\ell^1(\Ztd)$ into a non-commutative algebra with the delta-sequence $\delta$
as its unit element.

We tackle the problem to study the invertibility of twisted convolution
operators. Non-commutativity is the main subtle point
in this problem. In fact, the question when the mapping
\begin{equation*}
C_b:a\in\ell^1 \rightarrow \twc{a}{b} \in \ell^1
\end{equation*}
for some $b\in\ell^1$ is invertible and how we can compute the inverse
is more difficult than for a commutative setting. In particular, Wiener's
Lemma which deals with the problem that if, for some $b \in \ell^1$,
$C_b$ is invertible on $\ell^2$ then the inverse is generated from an element
again in $\ell^1$, has to be proven separately. An abstract and more general
proof of Wiener's Lemma for twisted convolution is given in \cite{GL03}.
Herein, we focus on a constructive method for studying the invertibility
of the twisted convolution with the rational parameter $q/p$.

In the following subsections we study the twisted convolution in a finite setting
and draw analogies for approaching the problem of invertibility of $C_b$ in the
general case.

\subsection{Twisted convolution on $\bf \Zp \times \Zp$}

%% As an example we describe the tw conv of the finite case.
%% This example also serves the purpose of building an analogy to
%% the more abstract part in which we reduce the tw conv
%% to a kind of "finite tw. conv".

In what follows we describe the twisted convolution on the finite
group $F = \Zp \times \Zp$. The standard (commutative) convolution
of two elements $f,g\in\C^{p\times p}$ is defined by
\begin{equation*}
(\conv{f}{g})_{m,n} =\sum_{k,l=0}^{p-1} f_{k,l}g_{m-k,n-l}\,,
\end{equation*}
where operations on indices is performed modulo $p$.

In analogy to the infinite case, we define the twisted convolution
$\twc{f}{g}$ of two elements $f,g\in\C^{p\times p}$ by
\begin{equation*}
   (\twc{f}{g})_{m,n} = \sum_{k,l = 0}^{p-1}  f_{k,l} g_{m-k,n-l}
   \om^{(m-k)l}
\end{equation*}
with $\om = e^{2\pi i q/p}$. For a fixed $g$, the twisted
convolution can be seen as a linear mapping $C_g:f \rightarrow \twc{f}{g}$
whose matrix $G$ is block circulant with $p$ blocks, i.e.,
\begin{equation*}
   G = C(G_0,G_{p-1},\dots,G_1)
   =\left (
   \begin{array}{cccc} G_0 & G_{p-1} & \cdots & G_1 \\
                            G_1 & G_0 & \cdots & G_2 \\
                \vdots & \vdots & & \vdots \\
                G_{p-1} & G_{p-2} & \cdots & G_0
   \end{array} \right )\,.
\end{equation*}
Each block has entries of the form
\begin{equation*}
  (G_j)_{kl} = \om^{jl}g_{j,k-l}\,.
\end{equation*}
Note that for the regular convolution each block is itself
circulant. For the invertibility of block circulant matrices we
apply a well known result from Fourier analysis.

\begin{lemma}{\rm \cite{Dav94}}
The matrix $G=C(G_0,G_{p-1},\dots,G_1)$ is invertible
if and only if every $\hat{G}_s = \sum_{r=0}^{p-1}e^{-2\pi i sr/p}G_r$,
$s=0,\dots,p-1$, is invertible. In this case
\begin{equation*}
  G^{-1} = C(H_0,H_{p-1},\dots,H_1)
\end{equation*}
where  $H_r = \frac{1}{p}\sum_{s=0}^{p-1}e^{2\pi i sr/p}(\hat{G}_s)^{-1}\,.$
\end{lemma}

By analyzing $\hat{G}_s$, we see that all blocks are unitary
equivalent, in the sense that
\begin{equation*}
   T_r\hat{G}_sT^*_r = \hat{G}_{s-qr}\,,
\end{equation*}
where $T_r$ denotes the unitary matrix with entries
\begin{equation*}
   (T_r)_{kl} =
   \left\{
   \begin{array}{cl}
     1 & \mbox{if} \quad p-r = l-k, \\
     0 & \mbox{else}\,.
   \end{array}
   \right.
\end{equation*}
Since $p$ and $q$ are relatively prime, we obtain all blocks by such a
unitary transformation. This implies that showing that if
$\hat{G}_0$ is invertible, then all $\hat{G}_s$ are invertible for
$s=1,\dots,p-1$. In other words, the $p\times p$ matrix $\hat{G}_0$
contains all the information about the invertibility of $C_g$.
An easy computation shows that the entries of $\hat{G}_0$
are given by
\begin{equation}\label{eq:entries}
   (\hat{G}_0)_{n,l} = \sum_{k=0}^{p-1}\om^{nl}g_{k,n-l}\,.
\end{equation}
We will later see that this observation motivates the matrix
algebra that we introduce to study the invertibility of the
twisted convolution.

Now, also all $\hat{G}^{-1}_s$ satisfy the same
unitary equivalence. It follows that we can read from $\hat{G}^{-1}_0$ the
element $g^{-1}$ which inverts the twisted convolution $f\rightarrow
\twc{f}{g}$, i.e., $\twc{g^{-1}}{g} = \twc{g}{g^{-1}} = \delta.$

The twisted convolution on $Z_p\times Z_p$ serves as analogy for modelling
the twisted convolution for the continuous and the finite dimensional
case.

%%%%%%%%%%%%%%%%%%%%%%%%%%%%%%%%%%%%%%%%%%%%%%%%%%%%%%%%%%%%%%%%%%%%%%%%

\section{Main results}

%% Here we describe the homomorphism between $\ell^1$ and a specific
%% matrix algebra in which we can explicitly prove Wiener's lemma.
%% Example: Finite dimensional case

Our aim is to find a way to describe those sequences that have an
inverse in $(\ell^1(\Ztd),\nat)$. To this end we divide the twisted
convolution into a finite sum of weighted normal convolutions of
sequences that have disjoint support. We define such a sequence
$a^{r,s}$ by
\begin{equation}\label{eq:coset-seq}
   (a^{r,s})_{k,l} = \left\{
   \begin{array}{ll}
      a_{k,l} \quad & \mbox{if} \quad (k,l) \equiv_p (r,s)\,, \\ [.2cm]
      0 & \mbox{else,}
   \end{array}
   \right.
\end{equation}
where $r,s \in \Zpd$. Obviously, $a^{r,s}$ is supported on the
coset $(r + p\Zd) \times (s + p\Zd)$ and $a = \sum_{r,s \in \Zpd}
a^{r,s}$. For a sequence $a$ having a coset support only for one
index, e.g., on $\Zd\times (s+p\Zd)$, we simply write
$a^{\cdot,s}$. We write $\equiv_p$ for denoting the equivalence of
integers modulo $p$. The idea of slitting a sequence into a sum of
sequences supported on cosets has first been introduced by
K.~Gr{\"o}chenig and W.~Kozek in \cite{GK97}.

\begin{lemma} \label{lem1} Let $a,b,c$ be in $\ell^1(\Ztd)$.
\begin{itemize}
\item[{\rm (a)}] For $r,s,u,v \in \Zpd$, $\conv{a^{r,s}}{b^{u,v}}$ is a
sequence supported on the coset $(u+r+p\Zd) \times (v+s+p\Zd)$.
\item[{\rm (b)}] If $c=c^{\cdot,0}$ is invertible in $(\ell^1,\ast)$,
then $c^{-1}$ is also supported on $\Zd \times p\Zd$.
\end{itemize}
\end{lemma}
\begin{proof}
Let $a^{r,s}, b^{u,v}$ be sequences in $\ell^1(\Ztd)$ and $k,l \in \Zpd$. Then
\begin{eqnarray*}
   (a^{r,s} \ast b^{u,v})_{k+p\Zd,l+p\Zd} &=& \sum_{m,n \in \Zd} (a^{r,s})_{m,n}
   (b^{u,v})_{k+p\Zd-m,l+p\Zd-n} \\
   &=& \sum_{m,n \in \Zpd} \sum_{(t,w) \equiv_p (m,n)} (a^{r,s})_{t,w}
   (b^{u,v})_{k+p\Zd-t,l+p\Zd-w}.
\end{eqnarray*}
Since $a^{r,s}$ is nonzero only for $(t,w) \equiv_p (r,s)$, and $b^{u,v}$
for $(k-t,l-w) \equiv_p (u,v)$, we obtain that $(k,l)$ has to be equivalent
to $(u+r,v+s)$ modulo $p$ for $a^{r,s} \ast b^{u,v}$ to be nonzero.

To show (b), let $c=c^{\cdot,0}$ be invertible and $e$ be its inverse. Then
\begin{equation*}
   \delta = c \ast e = c^{\cdot,0} \ast \big( \sum_{s \in \Zpd} e^{\cdot,s} \big )
   = \sum_{s \in \Zpd} c^{\cdot,0} \ast e^{\cdot,s},
\end{equation*}
where, by previous calculations, $c^{\cdot,0} \ast e^{\cdot,s}$ is a sequence
supported on $\Zd \times (s+p\Zd)$ for each $s \in \Zpd$. Since $\delta =
\sum_{s \in \Zpd} \delta^{\cdot,s}$, and elements of the sum have disjoint
supports, $c^{\cdot,0} \ast e^{\cdot,s} = \delta^{\cdot,s}$. But since
$\delta^{\cdot,s} = 0$ for $s \neq 0$ and $\delta^{\cdot,0} =
\delta$, we conclude that
\begin{equation*}
   c^{\cdot,0} \ast e^{\cdot,s} = \left \{
   \begin{array}{cc}
      \delta & s = 0\\ [.2cm]
      0 & s \neq 0
   \end{array} \right.
\end{equation*}
and therefore $e = e^{\cdot,0}$.
\end{proof}

With Definition (\ref{eq:coset-seq}), we obtain for $u,v \in \Zpd$
\begin{eqnarray*}
   (\twc{a}{b})_{u+p\Zd,v+p\Zd} & = &
   \sum_{k,l \in \Zd}a_{k,l}b_{u+p\Zd-k,v+p\Zd-l} \om^{\dotp{(u-k)}{l}} \\
   & = & \sum_{r,s \in \Zpd} \sum_{(k,l)\equiv_p (r,s)} a_{k,l}
   b_{u+p\Zd-k,v+p\Zd-l} \om^{\dotp{(u-r)}{s}} \\
   & = & \sum_{r,s \in \Zpd} \sum_{k,l \in \Zd}
   (a^{r,s})_{k,l}(b^{u-r,v-s})_{u+p\Zd-k,v+p\Zd-l} \om^{\dotp{(u-r)}{s}} \\
   & = & \sum_{r,s \in \Zpd} (a^{r,s}\ast
   b^{u-r,v-s})_{u+p\Zd,v+p\Zd} \om^{\dotp{(u-r)}{s}}\,.
\end{eqnarray*}
In a more compact notation we have
\begin{equation}\label{eq:comp}
   (\twc{a}{b})^{u,v} = \sum_{r,s \in \Zpd} a^{r,s} \ast b^{u-r,v-s}
   \om^{\dotp{(u-r)}{s}}\,.
\end{equation}
We observe now that the upper indices in (\ref{eq:comp}) behave like
a twisted convolution in $\Zpd \times \Zpd$. What changes is that we have
sequences as elements and standard convolution instead of multiplication.

Motivated by the block circulant structure of the twisted convolution on
$\Zp\times\Zp$ as described in the previous section, we introduce a new
matrix algebra which is isomorphic to $(\ell^1(\Ztd),\nat)$.

Before we do so, we fix an ordering of the elements from $\Zpd$.
Let $N= p^d$ and $\I=\{1,\ldots,N \}$. Then, to each $i \in \I$ we assign
an element $k_i$ from $\Zpd$ and set $k_1 = (0,\ldots,0)$. We will often write
$0$ instead of $k_1$.

Let $(\M,\cast)$ be an algebra of $p^d \times p^d$-matrices whose entries
are $\ell^1$-sequences and multiplication of two elements $A,B \in \M$ is
given by
\begin{equation*}
   (\mconv{A}{B})_{i,j} = \sum_{l \in \I} \conv{A_{i,l}}{B_{l,j}}
   \qquad i,j \in\I\,.
\end{equation*}
The identity element $\id$ is a matrix with $\delta$ sequences on the diagonal.

\begin{theorem}\label{thm1}
Let
\begin{equation*}
   \M_0 = \set{A \in \M}{A_{i,j} = \sum_{m \in \Zpd}
   \om^{\dotp{m}{k_j}} a^{m,k_i - k_j}, a \in \ell^1 ~\text{and}~ i,j \in \I}.
\end{equation*}
Then $\M_0$ is a subalgebra of $\M$.
\end{theorem}

\begin{proof}
Define a mapping $\phi \colon (\ell^1,\nat) \rightarrow (\M,\cast)$ by
\begin{equation}\label{eq:map-phi}
   (\phi(a))_{i,j} = \sum_{m \in \Zpd} \om^{\dotp{m}{k_j}} a^{m,k_i - k_j}.
\end{equation}
Then $\phi$ is linear, $(\phi(\delta))_{i,j} = \sum_{m \in \Zpd}
\om^{\dotp{m}{k_j}} \delta^{m,k_i-k_j} = \delta$ if $i=j$ and zero otherwise.
So $\phi(\delta) = \id$. We emphasize that the mapping $\phi$ has
been motivated by the matrix $\hat{G}_0$ described in the
previous section. For $i,j \in \I$,
\begin{eqnarray*}
  \big ( \phi(\twc{a}{b}) \big )_{i,j}
  &=& \sum_{m \in \Zpd} \om^{\dotp{m}{k_j}}
  (\twc{a}{b})^{m,k_i - k_j} = \sum_{m \in \Zpd} \om^{\dotp{m}{k_j}}
  \sum_{l,s \in \Zpd}
  \om^{(m-l)\cdot s} \conv{a^{l,s}}{b^{m-l,k_i- k_j-s}} \\
  &=& \sum_{m,l,s \in \Zpd} \om^{\dotp{m}{k_j+s})} \om^{\dotp{-l}{s}}
  \conv{a^{l,s}}{b^{m-l,k_i- k_j-s}} \\
  &=& \sum_{m,l,s \in \Zpd} \om^{\dotp{m}{s}} \om^{-l \cdot (s-k_j)}
  \conv{a^{l,s-k_j}}{b^{m-l,k_i-s}}\\
  &=& \sum_{m,l,s \in \Zpd} \om^{\dotp{l}{k_j}} \om^{\dotp{s}{(m-l)}}
  \conv{a^{l,s-k_j}}{b^{m-l,k_i-s}} \\
  &=& \sum_{s \in \Zpd} \conv{\big ( \sum_{l \in \Zpd}
  \om^{\dotp{l}{k_j}}a^{l,s-k_j} \big )}{\big ( \sum_{m \in \Zpd}
  \om^{\dotp{s}{(m-l)}} b^{m-l,k_i-s} \big )} \\
  &=& \sum_{n \in \I} \conv{\big ( \sum_{m \in \Zpd} \om^{\dotp{k_n}{m}}
  b^{m,k_i-k_n}\big)} {\big ( \sum_{l \in \Zpd} \om^{\dotp{l}{k_j}}
  a^{l,k_n-k_j} \big )} \\
  &=& \sum_{n \in \I} \conv{\phi(b)_{i,n}}{\phi(a)_{n,j}} \;= \;
  \big ( \mconv{\phi(b)}{\phi(a)} \big )_{i,j}\,.
\end{eqnarray*}
Therefore $\phi$ is an anti-homomorphism, that is,
\begin{equation*}
\phi(\twc{a}{b}) = \mconv{\phi(b)}{\phi(a)}.
\end{equation*}
Hence $\M_0$ is an algebra, being an image of an anti-homomorphism.
\end{proof}

Before stating the main theorem, we explore properties of elements of $\M_0$.
For $i,j \in \I$ and a matrix $A \in \M_0$ we define a new matrix $A(j,i)$
obtained from $A$ by substituting the $j$th row of $A$ with a vector of zeros
having $\delta$ on the $i$th position, and the $i$th column with a column of
zeros having $\delta$ on the $j$th position.

\begin{lemma} \label{lem2}
Let $A \in \mathcal{M}_0$. Then
\begin{itemize}
\item[\rm(a)] $\det (A)$ is a sequence supported on $\Zd \times p\Zd$.
\item[\rm(b)] $\det (A(1,i))$ is a sequence supported on
              $\Zd \times (k_i + p\Zd)$ for $i \in \I$.
\end{itemize}
\end{lemma}

\begin{proof}
Let $S_N$ be the group of permutations of the set $\I$. Then
\begin{eqnarray*}
  \det(A) &=& \sum_{\sigma \in S_N } (-1)^{\sigma} \prod_{i=1}^{N}
  A_{\sigma(i),i}= \sum_{\sigma \in S_N} (-1)^{\sigma} \prod_{i=1}^{N}
  \Big (\sum_{m_i \in \Zpd}\om^{k_i \cdot m_i} a^{m_i,k_{\sigma(i)}-k_i}\Big ) \\
  &=& \sum_{\sigma \in S_N} (-1)^{\sigma}
  \sum_{m_1,\ldots,m_{N} \in \Zpd} \om^{\sum_{i=1}^{N} \dotp{m_i}{k_i}}
  \underbrace{a^{m_1,k_{\sigma(1)}-k_1} \ast \cdots \ast
  a^{m_N,k_{\sigma(N)}-k_N}}_{G_{m_1,\ldots,m_N}}.
\end{eqnarray*}
Since $\sigma$ is a permutation of $\I$,
\begin{equation*}
   (k_{\sigma(1)}-k_1)+(k_{\sigma(2)}-k_2)+\cdots+(k_{\sigma(N)}-k_N)=0.
\end{equation*}
Therefore, by Lemma \ref{lem1}, $G_{m_1,\ldots,m_N}$ is a sequence
supported on the coset $(\sum_{i \in \I}m_i+p\Zd) \times p\Zd$.
Since $\sum_{i \in\I} m_i$ runs over all $\Zpd$, we see that $\det(A)$ is
supported on the coset $\Zd \times p\Zd$, i.e., $\det(A) = \det(A)^{\cdot,0}$.

In order to compute the support of $\det(A(1,i))$ for $i \in \I$,
let $S_{N-1}$ denote the group of permutations of $\{2,\ldots,N
\}$. Then for $i=1,\ldots,N$,
\begin{eqnarray*}
  \det(A(1,i))
  &=& (-1)^{i+1} \sum_{\sigma \in S_{N-1}}(-1)^{\sigma}
  A_{\sigma(2),1} \ast \cdots \ast A_{\sigma(i),i-1}
  \ast A_{\sigma(i+1),i+1} \ast \cdots \ast A_{\sigma(N),N}\\
  &=& (-1)^{i+1} \sum_{\sigma \in S_{N-1}} (-1)^{\sigma}
  \sum_{m_2,\ldots,m_N \in \Zpd} \om^{\dotp{m_2}{k_1} + \cdots +
  \dotp{m_i}{k_{i-1}}
  + \dotp{m_{i+1}}{k_{i+1}} + \cdots + \dotp{m_{N}}{k_N}} \times \\
  &\times& \underbrace{a^{m_2,k_{\sigma(2)}-k_1} \ast \cdots
  \ast a^{m_i,k_{\sigma(i)}-k_{i-1}}
  \ast a^{m_{i+1},k_{\sigma(i+1)}-k_{i+1}} \ast \cdots \ast
  a^{m_N,k_{\sigma(N)}-k_N}}_{G_{m_2,\ldots,m_N}}.
\end{eqnarray*}
Since $\sigma$ is a permutation of $\{2,\ldots,N\}$,
\begin{eqnarray*}
  &  (k_{\sigma(2)}-k_1) + \cdots + (k_{\sigma(i)}-k_{i-1}) +
  (k_{\sigma(i+1)}-k_{i+1}) + \cdots + (k_{\sigma(N)}-k_N) & \\
  &= \; (k_{\sigma(2)} + \cdots + k_{\sigma(N)}) - (k_1+k_2+\cdots+k_N) + k_i &\\
  &= \; (k_{\sigma(2)} + \cdots + k_{\sigma(N)}) - (k_2+\cdots+k_N) + k_i \; =
  \;k_i\,.&
\end{eqnarray*}
Therefore, by Lemma \ref{lem1}, $G_{m_2,\ldots,m_N}$ is supported on
$(\sum_{i=2}^{N} m_i + p\Zd) \times (k_i + p\Zd)$, and since each $m_i$
runs over all $\Zpd$, $\det(A(1,i))$ is supported on $\Zd \times (k_i+p\Zd)$.
That is, $\det(A(1,i)) = \det(A(1,i))^{\cdot,k_i}$.
\end{proof}

Now we are in the position to state and prove the main result

\begin{theorem}\label{th:wiener}{\rm [Wiener's Lemma for $\M_0$]}
Let $A \in\M_0$. If $A$ is invertible in $\M$, then $B=A^{-1} \in \M_0$.
\end{theorem}

\begin{proof}
Since $A \in \M_0$ is invertible, $\det(A)$ is an invertible sequence
in $(\ell^1,\ast)$, and there exists a matrix $\widetilde{B} \in \M$
such that $\mconv{A}{\widetilde{B}} = \id$. By Lemma \ref{lem2},
$\det(A) = \det(A)^{\cdot,0}$ and by Lemma \ref{lem1} its inverse,
$e = \det(A)^{-1}$, is also supported on the same coset, hence
$e = e^{\cdot,0}$. By Cramer's rule the inverse of $A$ is given by
\begin{equation*}
   \widetilde{B}_{i,j} = \conv{\det(A(j,i))}{e}.
\end{equation*}
We see that by Lemma \ref{lem2} (b), $\widetilde{B}_{i,1}$ is a
sequence supported on $\Zd \times (k_i+p\Zd)$. Let $b$ be a
sequence defined by
\begin{equation*}
   b = \widetilde{B}_{1,1}+\widetilde{B}_{2,1} + \ldots +
   \widetilde{B}_{N,1}.
\end{equation*}
Then $\widetilde{B}_{i,1} = \sum_{j \in \I} b^{k_j,k_i}$. Define
a new matrix, denoted by $B$, as
\begin{equation*}
   B_{i,j} = \sum_{m \in \Zpd} \om^{\dotp{m}{k_j}} b^{m,k_i-k_j}.
\end{equation*}
Then $B \in \M_0$ and we will show that $B=\widetilde{B}$, that is,
$B$ is the inverse of $A$.

Since $\widetilde{B}$ is the inverse of $A$,
\begin{eqnarray*}
   \id_{i,1} & = & \big ( \mconv{A}{\widetilde{B}} \big )_{i,1}
   = \sum_{j \in \I} \conv{A_{i,j}}{\widetilde{B}_{j,1}} \\ &=&
   \sum_{j \in \I} \sum_{m \in \Zpd} \om^{\dotp{m}{k_j}}
   \conv{a^{m,k_i-k_j}} {\widetilde{B}_{j,1}} \\
   &=& \sum_{j \in \I} \sum_{m \in \Zpd} \om^{\dotp{m}{k_j}} \conv{a^{m,k_i-k_j}}
   {\big ( \sum_{n \in \Zpd} b^{n,k_j} \big )} \\ &=& \sum_{j \in \I}
   \sum_{m \in \Zpd} \om^{\dotp{m}{k_j}} \sum_{n \in \Zpd}
   \conv{a^{m,k_i-k_j}}{b^{n,k_j}} \\
   &=& \sum_{m \in \Zpd}\sum_{j \in \I} \sum_{n \in \Zpd}
   \om^{\dotp{(m-n)}{k_j}} \conv{a^{m-n,k_i-k_j}}{b^{n,k_j}} \\ &=&
   \sum_{m \in \Zpd}G(m,k_i),
\end{eqnarray*}
where $G(m,k_i) = \sum_{j \in \I} \sum_{n \in \Zpd} \om^{\dotp{(m-n)}{k_j}}
 a^{m-n,k_i-k_j} \ast b^{n,k_j}$ is a sequence supported
on $(m+p\Z) \times (k_i+p\Z)$. Therefore, $G(k_1,k_1) = \delta$
and $G(m,k_i) = 0$ for $m \neq k_1$ and $i \neq 1$. Using the
above identity we will show that $\mconv{A}{B} = \id$, and by the
uniqueness of the inverse we will conclude that $B=\widetilde{B}$:
\begin{eqnarray*}
  & & \big ( \mconv{A}{B} \big )_{i,j} \; = \; \sum_{s \in \I} A_{i,s}
  \ast B_{s,j} \; = \\ &=& \sum_{s \in \I} \Big ( \sum_{m \in \Zpd}
  \om^{\dotp{m}{k_s}} a^{m,k_i-k_s} \Big) \ast \Big ( \sum_{n \in \Zpd}
  \om^{\dotp{n}{k_j}} b^{n,k_s-k_j} \Big ) \\
  &=& \sum_{m \in \Zpd} \sum_{s \in \I} \sum_{n \in \Zpd} \om^{\dotp{m}{k_s}}
  \om^{\dotp{n}{k_j}} a^{m,k_i-k_s} \ast b^{n,k_s-k_j}\\
  &=& \sum_{m \in \Zpd} \sum_{s \in \I} \sum_{n \in \Zpd} \om^{\dotp{(k_s+k_j)}{m}}
  \om^{\dotp{n}{k_j}} a^{m,k_i-k_j-k_s} \ast b^{n,k_s}\\ &=&
  \sum_{m \in \Zpd} \sum_{s \in \I} \sum_{n \in \Zpd} \om^{\dotp{(k_s+k_j)}{(m-n)}}
  \om^{\dotp{n}{k_j}} a^{m-n,(k_i-k_j)-k_s} \ast b^{n,k_s} \\
  &=& \sum_{m \in \Zpd} \om^{\dotp{m}{k_j}} \sum_{s \in \I} \sum_{n \in \Zpd}
  \om^{\dotp{(m-n)}{k_s}} a^{m-n,(k_i-k_j)-k_s} \ast b^{n,k_s} \\
  &=& \sum_{m \in \Zpd} \om^{\dotp{m}{k_j}} G(m,k_i-k_j) \; = \;
  \left \{ \begin{array}{cc}
     \delta & k_i - k_j = k_1 ~ \Leftrightarrow ~ i=j;\\ [.2cm]
     0 & k_i - k_j \neq k_1 ~ \Leftrightarrow ~ i \neq j;
  \end{array} \right.
\end{eqnarray*}
Hence,  $\mconv{A}{B} =I$.
\end{proof}

Theorem \ref{th:wiener} provides the key result to study
invertibility of twisted convolution. Indeed, for a given sequence
$a$ in $\ell^1$ we look at the corresponding matrix $A = \phi(a)$
as defined in (\ref{eq:map-phi}). If $A$ is invertible in
$(\M,\cast)$, which can be checked showing that the determinant is
invertible in $(\ell^1,\ast)$, then its inverse $A^{-1}$ is of the form
$\phi(b)$ for another element $b$ in $\ell^1$. This element $b$, in turn,
provides the inverse of $a$ in $(\ell^1,\nat)$.

The approach is constructive in the sense that algebraic methods such as
Cramer's Rule can be applied to find the inverse of $A$. Then, the
sequence $b$ can simply be read from the entries of $A^{-1}$ according
to the mapping $\phi$. In particular for small $p$ and $d$ this method
leads to fast inversion schemes for the twisted convolution operator.
In the last section we will show explicitly how this works in the case
of $d=1$.

%%%%%%%%%%%%%%%%%%%%%%%%%%%%%%%%%%%%%%%%%%%%%%%%%%%%%%%%%%%%%%%%%%%%%%%%

\section{Twisted convolution and Gabor analysis}

%% A motivation for the twisted convolution !!

Central objects in time frequency analysis are
modulation and translation operators. Although most of the upcoming
notation can be given in the more general setting of locally compact
Abelian groups we restrict ourselves to $\Rd$ in order to simplify
the readability of this article.

For $x,\om\in\Rd$ we define the translation operator and
the modulation operator on $L^2(\Rd)$ by
\begin{eqnarray*}
T_x f(\cdot) &=& f(\cdot - x)\,, \\
M_\om f(\cdot) &=& e^{2\pi i \om \cdot} f(\cdot)\,,
\end{eqnarray*}
respectively. Many technical details in time-frequency analysis are
linked to the commutation law of the translation and modulation operator,
namely,
\begin{equation}\label{eq:commlaw}
 M_\om T_x\, = \, e^{2\pi i \dotp{x}{\om}}T_x M_\om  \,.
\end{equation}
The time-frequency shift for $x,\om\in\Rd$ is denoted by
\begin{equation*}
 \pi(x,\om) = T_x M_\om.
\end{equation*}
It follows from (\ref{eq:commlaw}) that
\begin{equation}\label{eq:commlaw2}
\pi(x_1,\om_1)\pi(x_2,\om_2) \, = \, e^{2\pi i \dotp{x_2}{\om_1}}
\pi(x_1+x_2,\om_1+\om_2)\,.
\end{equation}
This shows that time-frequency shifts almost allow a group
structure. Incorporating the additional phase factor into a more
extended group law leads to the so-called Heissenberg group. For
more details about this topic, the reader is referred to
\cite{F89}.

Gabor analysis deals with the problem of decomposing and reconstructing
signals according to a special basis system which consists of regular
time-frequency shifts of a single so-called window function \cite{FS98, FS03}.
Let $\La$ be a time-frequency lattice, i.e., a discrete subgroup of the
time-frequency plane $\Rtd$, and let $g$ be in $L^2(\Rd)$. Then we
define a Gabor system $\G(g,\La)$ by
\begin{equation*}
   \G(g,\La) \, = \, \set{\pi(\la)g}{\la\in\La}\,.
\end{equation*}
We associate with this Gabor system the positive operator
\begin{equation*}
   S:f\in L^2 \rightarrow Sf = \sum_{\la\in\La}
   \inner{f}{\pi(\la)g}\pi(\la)g \,.
\end{equation*}
If the operator $S$ is bounded and invertible on $L^2(\Rd)$,
then $\G(g,\La)$ is called a frame and $S$ the associated frame
operator, cf.~\cite{CH03}.

Many studies in Gabor analysis are devoted to the frame operator
\cite{Gro01}. In what follows we will describe the so-called
Janssen representation of such operators. To this end we need
the notion of the adjoint lattice, i.e.,
\begin{equation*}
   \aLa \, = \, \set{\ala\in\Rtd}{\pi(\la)\pi(\ala)=
   \pi(\ala)\pi(\la),\, \la\in\La}\,.
\end{equation*}
In \cite{DLL95,FK98,Jan95} it is shown that the frame operator $S$
satisfies Janssen representation,
\begin{equation} \label{eq:jans-repr}
   S \, = \, \sum_{\ala\in\aLa} \inner{g}{\pi(\ala)g}\pi(\ala)\,.
\end{equation}
At this point, the question arises if we can deduce the invertibility
of the operator $S$ from the Janssen coefficients
$(\inner{g}{\pi(\ala)g})$. It is known from frame theory that if $S$
is invertible, then its inverse is of the same type, that is,
it also has a Janssen representation.

In order to better understand the main ingredients of this problem
we transfer the model to an operator algebra. To this end
we restrict our discussion to so-called separable lattices of
the form
\begin{equation*}
       \La \, = \, \al \Zd \times \be \Zd
\end{equation*}
for some fixed positive numbers $\al$ and $\be$. An easy computation
based on (\ref{eq:commlaw2}) shows that
\begin{equation*}
       \aLa \, = \, \be^{-1} \Zd \times \al^{-1}  \Zd \,.
\end{equation*}

We define the operator algebra $\A$ as in \cite{GL03} by
\begin{equation*}
       \A \, = \, \set{S = \sum_{k,l\in\Zd}a_{k,l}
       \pi(\be^{-1}k,\al^{-1}l)}{a=(a_{k,l})\in\ell^1(\Ztd)}\,.
\end{equation*}
The restriction to $\ell^1$-sequences guarantees absolute convergence
of the sum of time-frequency shifts. Let $\ka$ be the mapping
\begin{equation*}
    \ka:a\in\ell^1 \rightarrow \ka(a) =  \sum_{k,l\in\Zd}a_{k,l}
       \pi(\be^{-1}k,\al^{-1}l)\in\A \,.
\end{equation*}
Then, as already observed in \cite{Jan95}, we have
\begin{equation*}
    \ka(a)\ka(b) \, = \, \ka(\twc{a}{b})
\end{equation*}
and $\ka(\delta) = \id$ where $\delta$ and $\id$ denote the Dirac
sequence and the identity operator, respectively. Both represent
the unit element of the corresponding algebra. It follows that $\ka$
is an algebra homomorphism from $(\ell^1(\Ztd),\nat)$ to $\A$, and
invertibility of an element in $\A$ can be transferred to the
invertibility of the associated $\ell^1$-sequence with respect to
the twisted convolution.

It is important to observe, that all the results go through also for
weighted $\ell^1$-spaces. These facts are used to design dual
Gabor windows of a special type, cf.~\cite{GL03}.

In the following section we give an example of how this approach can be
explicitly used in Gabor analysis of one-dimensional signals.

\noindent{\bf Remark.} The above results, with the help of metaplectic
operators, carry over to the more general class of lattices,
called symplectic lattices. A lattice $\La_s \subseteq \Rtd$ is
called symplectic, if one can write $\La_s = \D \La$ where $\La$
is a separable lattice and $\D \in GL_{2d}(\R)$. To every $\D \in
GL_{2d}(\R)$, there corresponds a unitary operator $\mu(\D)$,
called metaplectic, acting on $L^2(\Rd)$. One can show that a
Gabor system on a symplectic lattice is unitary equivalent to a
Gabor system on a separable lattice under $\mu(\D)$, and
\begin{equation*}
S_{g}^{\La_s} = \mu(\D)^{-1} S_{\mu(\D)g}^{\La} \mu(\D).
\end{equation*}
Hence, to analyze the invertibility of a frame operator $S$
associated to the window function $g \in L^2(\Rd)$ and symplectic
lattice $\La_s$, it suffices to analyze a frame operator
associated to $\mu(\D)g$ and a separable lattice $\La$.
For more details see \cite{Gro01}.

%%%%%%%%%%%%%%%%%%%%%%%%%%%%%%%%%%%%%%%%%%%%%%%%%%%%%%%%%%%%%%%%%%%%%%%%

\section{Application to one-dimensional signal space}

%% Herein we shortly describe the 1D algorithm that is
%% studied in the sp-paper in more details.

In this section, we briefly describe how the presented inversion scheme
applies to Gabor frame operators in a one-dimensional setting.
A more detailed discussion also for finite dimensional signals
is described in \cite{MWE05}.

Assume $d=1$. Let $a$ be in $\ell^1(\Z^2)$ and $\al, \be$ be constants such
that $\al\be = p/q$ with $p,q$ relative prime. Set
\begin{equation*}
 \ka(a) = \sum_{k,l\in\Z}a_{k,l}\pi(\be^{-1}k,\al^{-1}l)\,.
\end{equation*}
In order to verify if $\ka(a)$ is invertible on $L^2(\R)$ we simply look
at the coefficient sequence $a$ and check whether $a$ is invertible in
$(\ell^1(\Z^2),\nat)$. To this end, we apply the above results and
switch to the matrix $A$ whose entries are defined by
\begin{equation*}
  A_{i,j} = \sum_{m=0}^{p-1} \om^{mj}a^{m,i-j}\,,
\end{equation*}
with $\om = e^{2\pi i q/p}$. Next, we need to show that the matrix $A$
is invertible in $(\M,\cast)$. For example, we can calculate the
determinate which is a sequence in $\ell^1$ and show that it is
invertible in $(\ell^1,\ast)$.

Assume that the determinant of $A$ is invertible. We denote its inverse
by $e$. By Cramer's Rule, we compute the first column of the inverse
matrix $B$ of $A$ as
\begin{equation*}
  B_{k,0} = \conv{\det A(0,k)}{e}\,,
\end{equation*}
for $k=0,\dots,p-1$. Then
\begin{equation*}
  b = \sum_{k=0}^{p-1}B_{k,0}
\end{equation*}
provides the inverse sequence of $a$ which, in turns, gives
$\ka(a)^{-1} = \ka(b)$.

Note that for $p=1$, the twisted convolution turn into normal convolution
and we can simply apply the standard Fourier inversion scheme of
sequences in $(\ell^1(\Z^2),\ast)$ since in this case the matrix $A$ reduces
to the sequence $a$.

%%%%%%%%%%%%%%%%%%%%%%%%%%%%%%%%%%%%%%%%%%%%%%%%%%%%%%%%%%%%%%%%%%%%%%%%

\section*{Acknowledgements} The authors gratefully acknowledge support from the
Ollendorff Minerva Center and from the European Union's Human Potential Programme,
under the contract HPRN-CT-2003-00285 (HASSIP). We would also like to thank
Karlheinz Gr{\"o}chenig and Yehoshua Y.~Zeevi for many fruitful discussions.

%%%%%%%%%%%%%%%%%%%%%%%%%%%%%%%%%%%%%%%%%%%%%%%%%%%%%%%%%%%%%%%%%%%%%%%%
%%%%%%%%%%%%%%%%%%%%%%%%%%%%%%%%%%%%%%%%%%%%%%%%%%%%%%%%%%%%%%%%%%%%%%%%

\bibliography{twc}
\bibliographystyle{plain}
%\bibliographystyle{alpha}

%%%%%%%%%%%%%%%%%%%%%%%%%%%%%%%%%%%%%%%%%%%%%%%%%%%%%%%%%%%%%%%%%%%%%%%%
%%%%%%%%%%%%%%%%%%%%%%%%%%%%%%%%%%%%%%%%%%%%%%%%%%%%%%%%%%%%%%%%%%%%%%%%

\end{document}